\documentclass[11pt]{article}
\usepackage{amsmath,amssymb,graphics,amsthm}

\newtheorem{lemma}{Lemma}

\newtheorem{proposition}{Proposition}
\newtheorem{question}{Question}
\newtheorem{corollary}{Corollary}
\newtheorem{conjecture}{Conjecture}
\newtheorem{theorem}{Theorem}
\newtheorem{definition}{Definition}

\def\proof{{\bf Proof. }}
\begin{document}
\title{{\sc On Nash's 4-sphere and Property 2R}}
\author{Motoo Tange}
\date{}
\maketitle
\abstract
D.Nash defined a family of homotopy 4-spheres in \cite{[n]}.
Proving that his manifolds ${\mathcal S}_{m,n,m',n'}$ are all real $S^4$,
we find that they have handle decomposition with no 1-handles, two 2-handles and two 3-handles.
The handle structures give new potential counterexamples of Property 2R conjecture.
\section{Introduction}
The smooth Poincar\'e conjecture in 4-dimension is still open.
Though many people \cite{[ak],[cs]} have proposed potential counterexamples, what some 
of them are standard $S^4$ was proved \cite{[a1],[g1],[g2],[g3]}.
D. Nash \cite{[n]} also proposed potential counterexamples of the conjecture.
Most recently S.Akbulut \cite{[a2]} proved that the manifolds are all standard.
In the article I will also give an alternative proof and furthermore remark some handle decompositions
appeared there.

Nash's manifolds are constructed by log transformations along four tori in some 4-manifold.
Hence we will give a brief review of the sugery.
For the remark of the handle decomoposition as stated above we introduce notions; Property nR, generalized
Property R.

\subsection{Log transformation.}
Here we review the notation of the log transformation.
Let $T\subset X^4$ be a torus embedding with the trivial normal bundle $\nu(T)=D^2\times T$ in 4-manifold $X$.
Removing the neighborhood, we reglue it with the map $\varphi:\partial D^2\times T^2\to \partial \nu(T)$ satisfying
$$\varphi(\partial D^2\times \{\text{pt}\})=p\mu+q\gamma,$$
where $\mu$ is the meridian of $T$ and $[\mu]$ is a primitive element of $H_1(T)$, 
so that we obtain a manifold.
\begin{definition}
The surgery whose gluing map is $\varphi$ as above 
$$X-\nu(T)\cup_\varphi(D^2\times T)$$
is called the (p/q)-log transformation along $T$ with direction $\gamma$.
\end{definition}
\subsection{Generalized Property R Conjecture.}
Property R conjecture was proved by Gabai \cite{[ga]}.
M. Scharlemann and A. Thompson in \cite{[st]} generalized  Property R as follows.
\begin{definition}[\cite{[st]}]
We say that a knot $K$ has Property nR if $K$ satisfies the following property.
If any $n$-component link $L$ containing $K$ as a component
yields $\#^nS^1\times S^2$ by an integral Dehn surgery, then after some handle slidings
the framed link can be reduced to the $n$-component unlink.
\end{definition}
The case where $n=1$ is equivalent to original Property R.
\begin{conjecture}[Generalized Property R Conjecture]
All knot admit Property nR for any $n\ge 1$.
\end{conjecture}
The generalized Property R conjecture is still open.
The homotopy 4-spheres by D.Nash in \cite{[n]} are standard, however we show that
diagrams coming from handle decompositions might be counterexamples of the generalized Property R conjecture.

\begin{figure}[htpb]
\begin{center}
\includegraphics{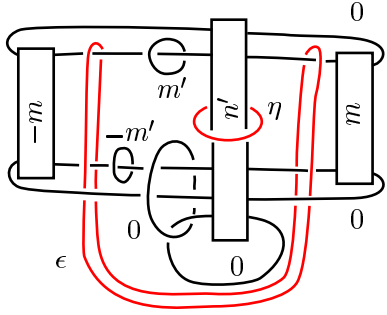}
\caption{Examples which might not have Property 2R.}
\label{pcp2}
\end{center}\end{figure}

We can find Figure~\ref{pcp2} along the way of proof that Nash's manifolds are standard (Theorem~\ref{main}).
The framed link with black color is a presentation of $S^3$.
Each box stands for the full twist by the number in the box.
The $0$-surgery along the 2-component link (red color) gives rise to $\#^2S^2\times S^1$, because
the framed link gives a handle diagram of Nash's homotopy 4-sphere, which are indeed real $S^4$ (Corollary~\ref{cor1}).
\begin{question}
Are $\eta$ and $\epsilon$ in Figure~\ref{pcp2} examples not having Property 2R for any non-zero integers $m,n,m'$?
\end{question}

\section*{Acknowledgements}
The whole work of the paper was done during my visiting of The University of Texas at Austin in 2011.
I deeply appreciate Professor Cameron Gordon's acceptance of my visiting, his, Professor John Luecke's and Dr. Cagri Karakurt's hospitality,
many pieces of chalk and wide blackboard to draw some pictures.
Also, this research was supported by Grant-in-Aid for JSPS Fellow(21-1458).
\section{Nash's manifolds.}
\begin{figure}[thbp]
\begin{minipage}{.5\textwidth}
\begin{center}
\includegraphics{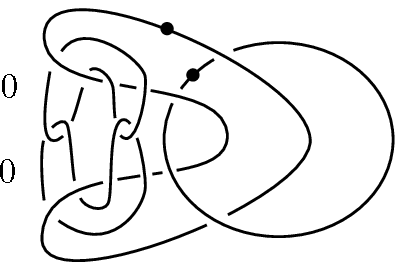}
\caption{$A$.}
\label{A}
\end{center}
\end{minipage}
\begin{minipage}{.5\textwidth}
\begin{center}
\includegraphics{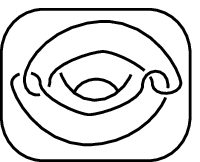}
\caption{Bing tori.}
\label{B}
\end{center}
\end{minipage}
\end{figure}

D. Nash in \cite{[n]} defined a new family of homotopy 4-spheres as follows.
Let $A$ be a 4-manifold with the handle diagram Figure~\ref{A}.

Since $A$, as in \cite{[fs]}, is constructed by attaching two 2-handles to 
$D^2\times T^2_0\subset (D^2\times S^1)\times S^1$ along Bing tori of $D^2\times S^1$,
$A$ is also included a Bing tori $B_T$ as in Figure~\ref{B},
where $T_0^2$ is the punctured torus.
As a fundamental fact $0$-surgery along $B_T$ yields $T^2_0\times T^2_0$.
Here $0$-surgery means $(0/1)$-log transformation.
In addition two tori reglued by the surgery are $T_1=S^1_\alpha\times S^1_\gamma$ and $T_2=S^1_\beta\times S^1_\gamma$,
where $S^1_\alpha,S^1_\beta,S^1_\gamma$ are generating circles of $T^2_0\subset T^2$.
In other words tori (0/1)-surgery along $T^2_1$ and $T^2_2$ yield $A$ back.

Now we take two copies of $T^2_0\times T^2_0$ to glue the boundaries $S^1\times T^2_0\cup T^2_0\times S^1$ each other by 
gluing map $\phi:S^1\times T^2_0\cup T^2_0\times S^1\to S^1\times T^2_0\cup T^2_0\times S^1$ such that 
the two components are exchanged.
We call the resulting manifold $X$.
Such a construction is gereralized to $n$-component case by Fintushel-Stern in \cite{[fs]}, and
it is called {\it pinwheel construction}.
Here we define the $(m/1)$-surgery of $T_0^2\times T_0^2$ along $T_1$ with direction $S^1_\alpha$ and simultaneously
the $(n/1)$-surgery along $T_2$ with direction $S^1_\beta$ to be $X_{m,n}$.
We define, by the same gluing map $\phi$, $X_{m,n}\cup_\phi X_{m',n'}$ to be ${\mathcal S}_{m,n,m',n'}$.
From the construction immediately we have the following.
\begin{lemma}
\label{exchange}
For any integer $m,n,m',n'$ we have the following diffeomorphism
$${\mathcal S}_{m,n,m',n'}\cong {\mathcal S}_{m',n',m,n}.$$
\end{lemma}

Here Nash got a result (Theorem~3.2 in \cite{[n]}), which
the manifolds ${\mathcal S}_{m,n,m',n'}$ are all homotopy 4-spheres.
Namely ${\mathcal S}_{m,n,m',n'}$ are candidates of the 4-dimensional smooth Poincar\'e conjecture.
Are the manifolds diffeomorphic to standard $S^4$?
Here we give an affirmative answer for the question.
\begin{theorem}[Nash's manifolds are standard.]
\label{main}
The manifolds ${\mathcal S}_{m,n,m',n'}$ are all diffeomorphic to the standard 4-sphere.
\end{theorem}
S.Akbulut independently proved the same result in \cite{[a2]}.

\section{Handle decomposition of ${\mathcal S}_{m,n,m',n'}$.}
\subsection{The diagram of $X_{m,n}$.}
\begin{lemma}
A handle decommposition of $X_{m,n}$ is Figure~\ref{upd}.
\end{lemma}

\begin{figure}[thbp]
\begin{center}
\includegraphics{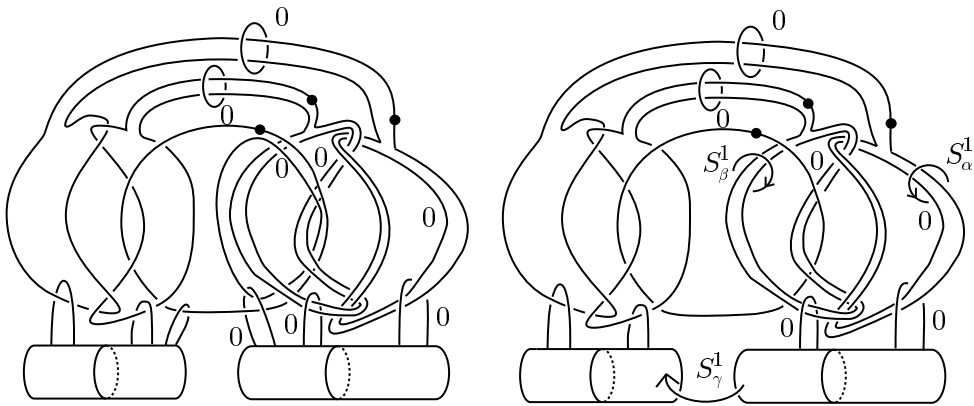}
\caption{$T^4$ and $T^2_0\times T^2_0$.}
\label{T4-T0T0}
\end{center}
\end{figure}

\proof
First the picture of $T^4$ is the left of Figure~\ref{T4-T0T0}.
Recall that $T^3$ is obtained from $0$-surgery along the Borromean ring.
Since $T^2_0\times T^2_0$ is obtained by removing $D^2_0\times T^2\subset \ T^2\times D^2$
and $T_0^2\times D^2_0\subset T^2_0\times T^2$, the diagram is the right of Figure~\ref{T4-T0T0}.
Since $(m/1)$ and $(n/1)$-log transformation correspond to the $(0,1/m,1/n)$ surgery over the Borromean ring,
we get Figure~\ref{upd} as a diagram of $X_{m,n}$.
\qed

The 3 directions of the right in Figure~\ref{T4-T0T0} represent $S^1_\alpha,S^1_\beta,S^1_\gamma$ generating circles above.
\subsection{Upside-down of $X_{m',n'}$.}
\begin{figure}[thbp]
\begin{center}
\includegraphics{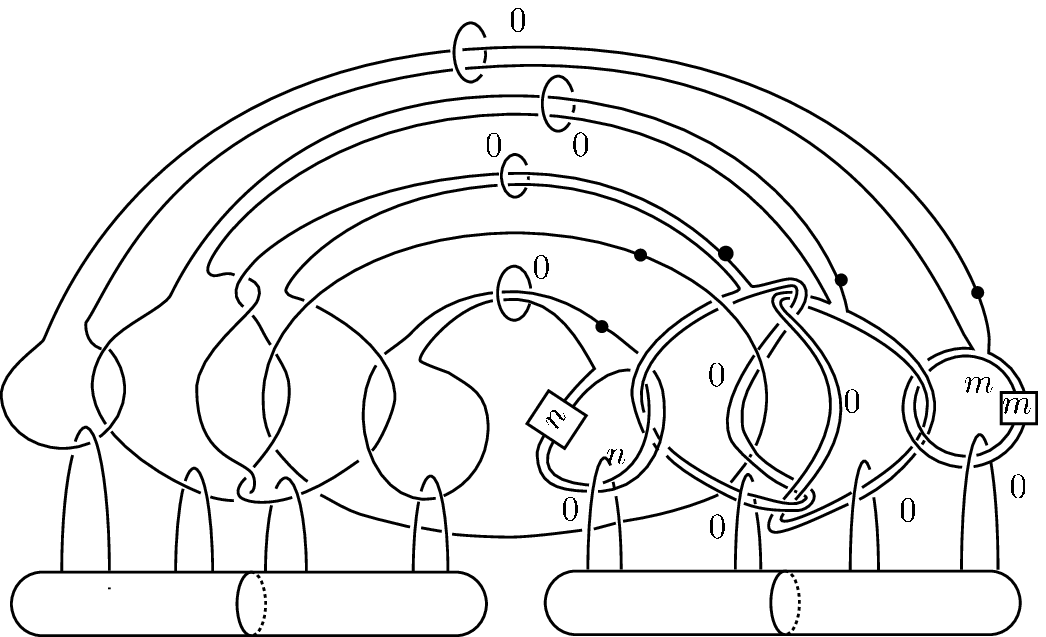}
\caption{$X_{m,n}$}
\label{upd}
\end{center}
\end{figure}

\begin{figure}[thbp]
\begin{center}
\includegraphics{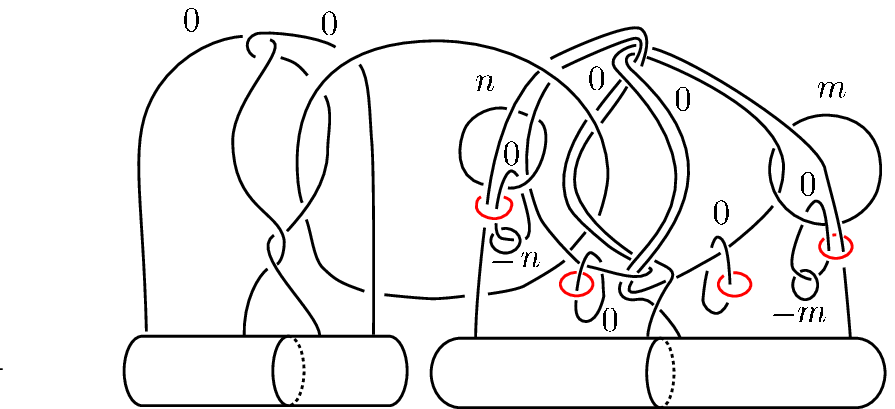}
\caption{}
\label{upd6}
\end{center}\end{figure}

Next we perform the upside down of the manifolds $X_{m',n'}$.
The right four 2-handles in Figure~\ref{upd} which are along two components of Borromean ring and the two meridiands
are, as each runs throught the adjacent 1-handles once, canceled each other.
In addition the top four 2-handles are isotopic to trivial unknots on the boundary and they are canceled out with
four 3-handles.
Then attaching dual 2-handles are the meridian for the bottom four 2-handles as in Figure~\ref{upd}.
Here we present the dual 2-handle by red lines.
Then by handle sliding we get the diagram Figure~\ref{upd6}.

\begin{figure}[thbp]
\begin{center}
\includegraphics{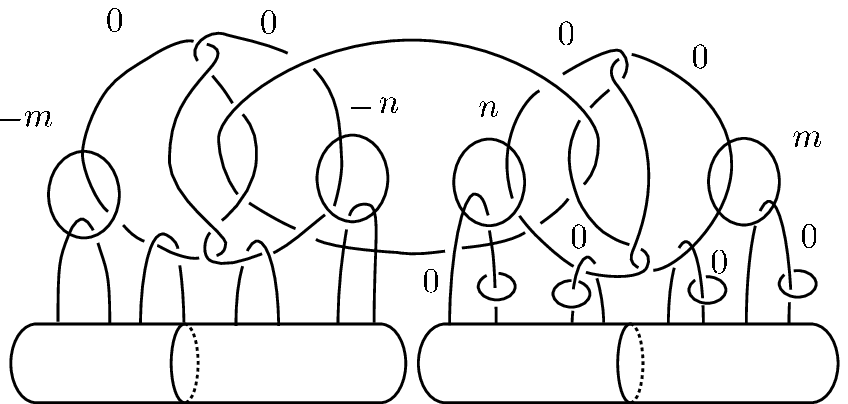}
\caption{}
\label{upd2}
\end{center}
\end{figure}

\begin{figure}[thbp]
\begin{center}
\includegraphics{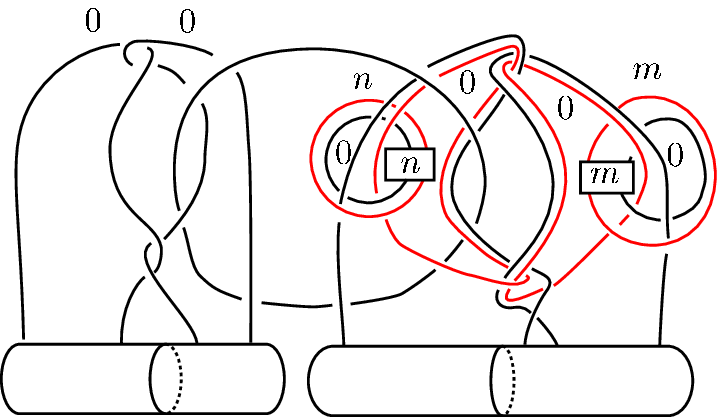}
\caption{}
\label{upd7}
\end{center}\end{figure}

\begin{figure}[thbp]
\begin{center}
\includegraphics{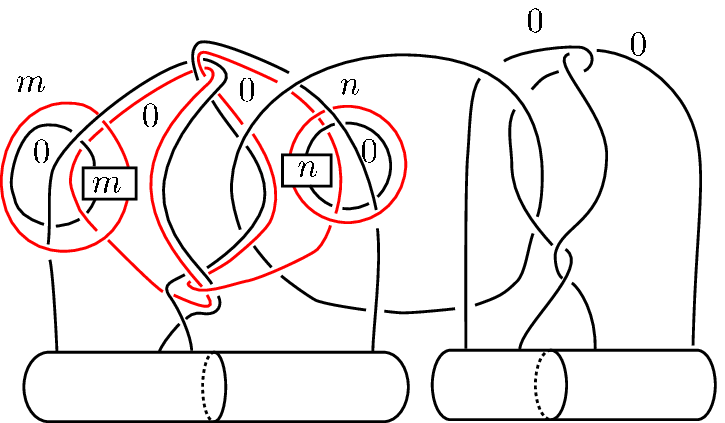}
\caption{}
\label{upd8}
\end{center}\end{figure}

\begin{figure}[thbp]
\begin{minipage}{.5\textwidth}
\begin{center}
\includegraphics{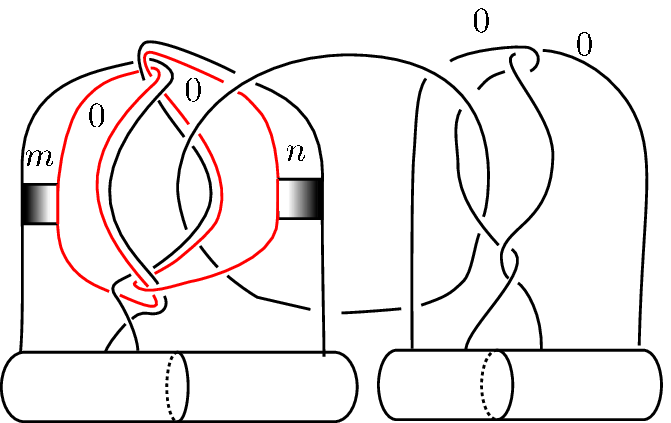}
\caption{}
\label{upd10}
\end{center}
\end{minipage}
\begin{minipage}{.5\textwidth}
\begin{center}
\includegraphics{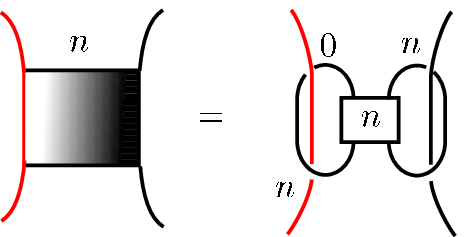}
\caption{}
\label{abb}
\end{center}
\end{minipage}
\end{figure}

\begin{figure}[thbp]
\begin{minipage}{.5\textwidth}
\begin{center}
\includegraphics{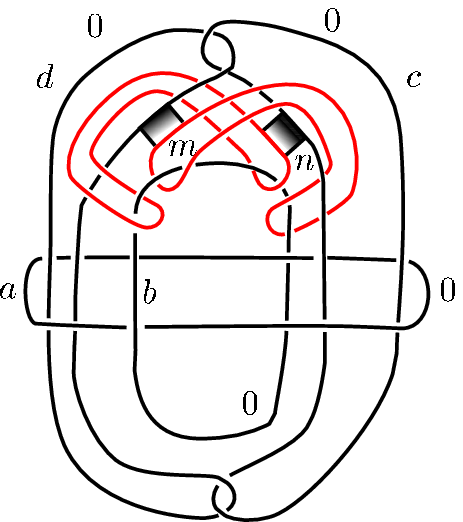}
\caption{}
\label{upd16}
\end{center}
\end{minipage}
\begin{minipage}{.5\textwidth}
\begin{center}
\includegraphics{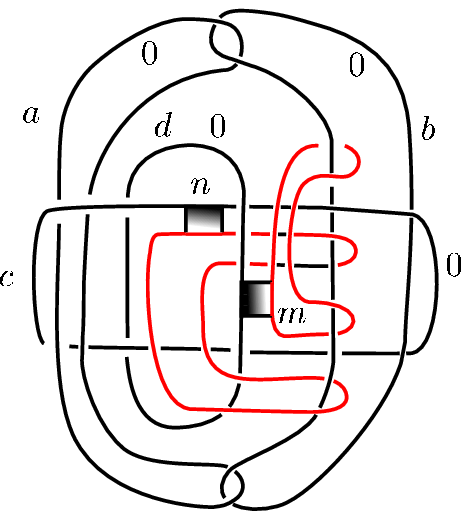}
\caption{}
\label{upd14}
\end{center}
\end{minipage}
\end{figure}

\begin{figure}[htpb]
\begin{center}
\includegraphics{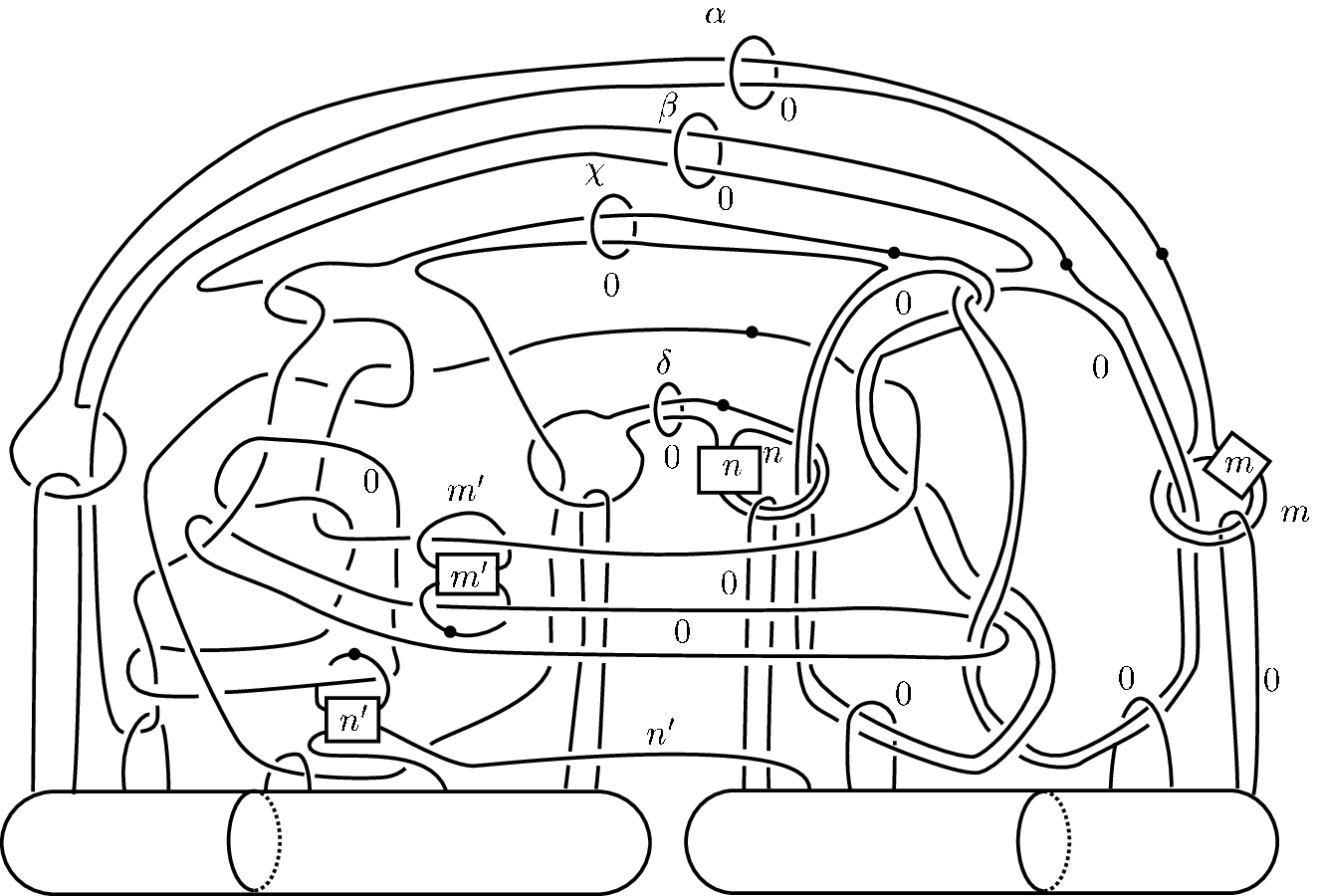}
\caption{${\mathcal S}_{m,n,m',n'}$.}
\label{upd22}
\end{center}\end{figure}

In addition several handle slides give Figure~\ref{upd7} and \ref{upd8}.
Here replacing the two handles as in Figure~\ref{abb} we get Figure~\ref{upd10}.
Using the notation and isotopy we get Figure~\ref{upd16} and keep track of the red two handles 
by the symmetry that exchanges the pair of link $(a,b)$ to $(c,d)$,
hence we get Figure~\ref{upd14}.
Keeping track of the diagram by the converse motion (Figure~\ref{upd16}-\ref{upd10}-\ref{upd8}-\ref{upd7}-\ref{upd6}-\ref{upd2}-\ref{upd})
from the diagram in the form, we get Figure~\ref{upd22}.

\section{Handle calculus of ${\mathcal S}_{m,n,m',n'}$.}
\begin{figure}[htpb]
\begin{center}
\includegraphics{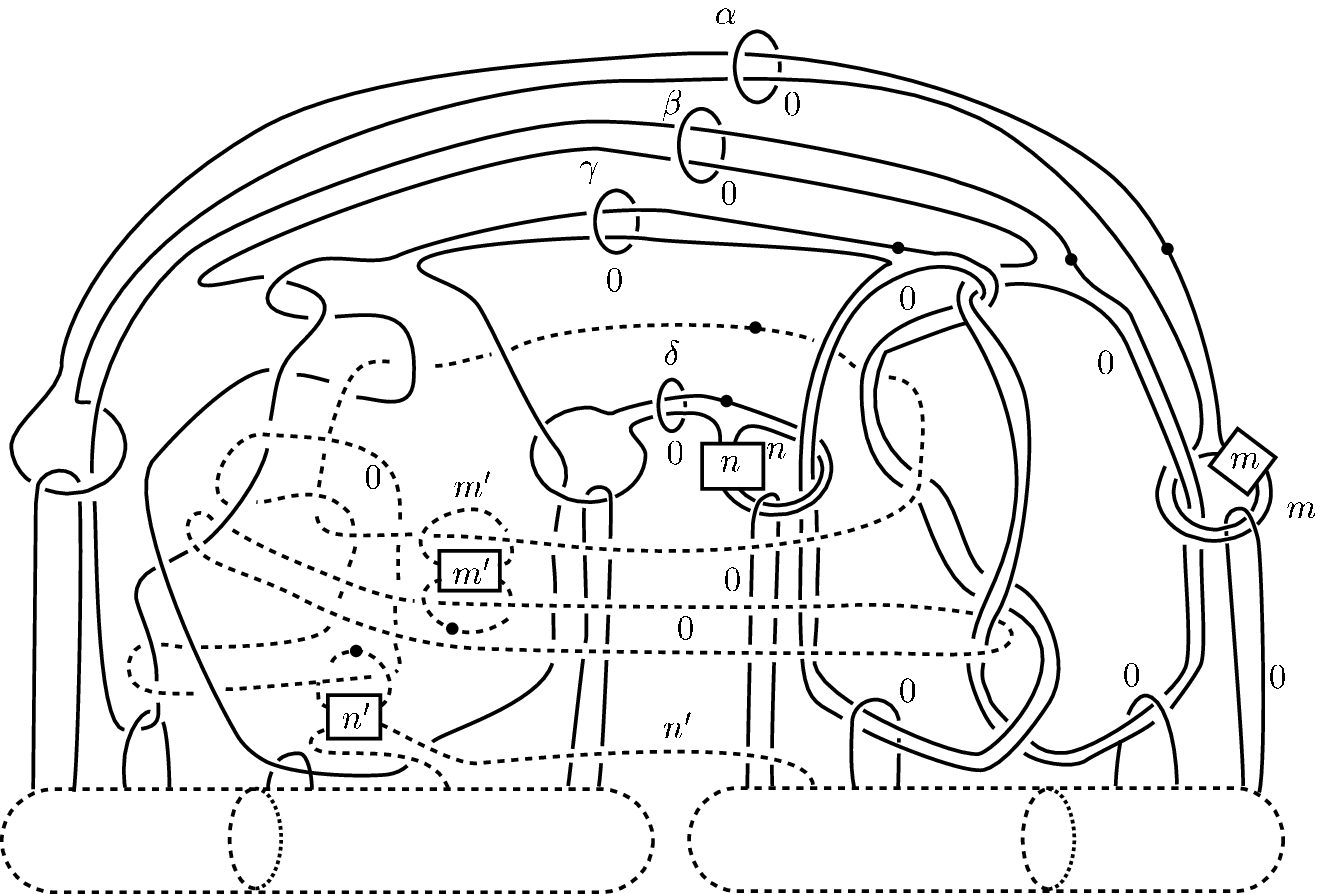}
\caption{${\mathcal S}_{m,n,m',n'}$ canceled.}
\label{can1}
\end{center}\end{figure}

\begin{figure}[htpb]
\begin{center}
\includegraphics{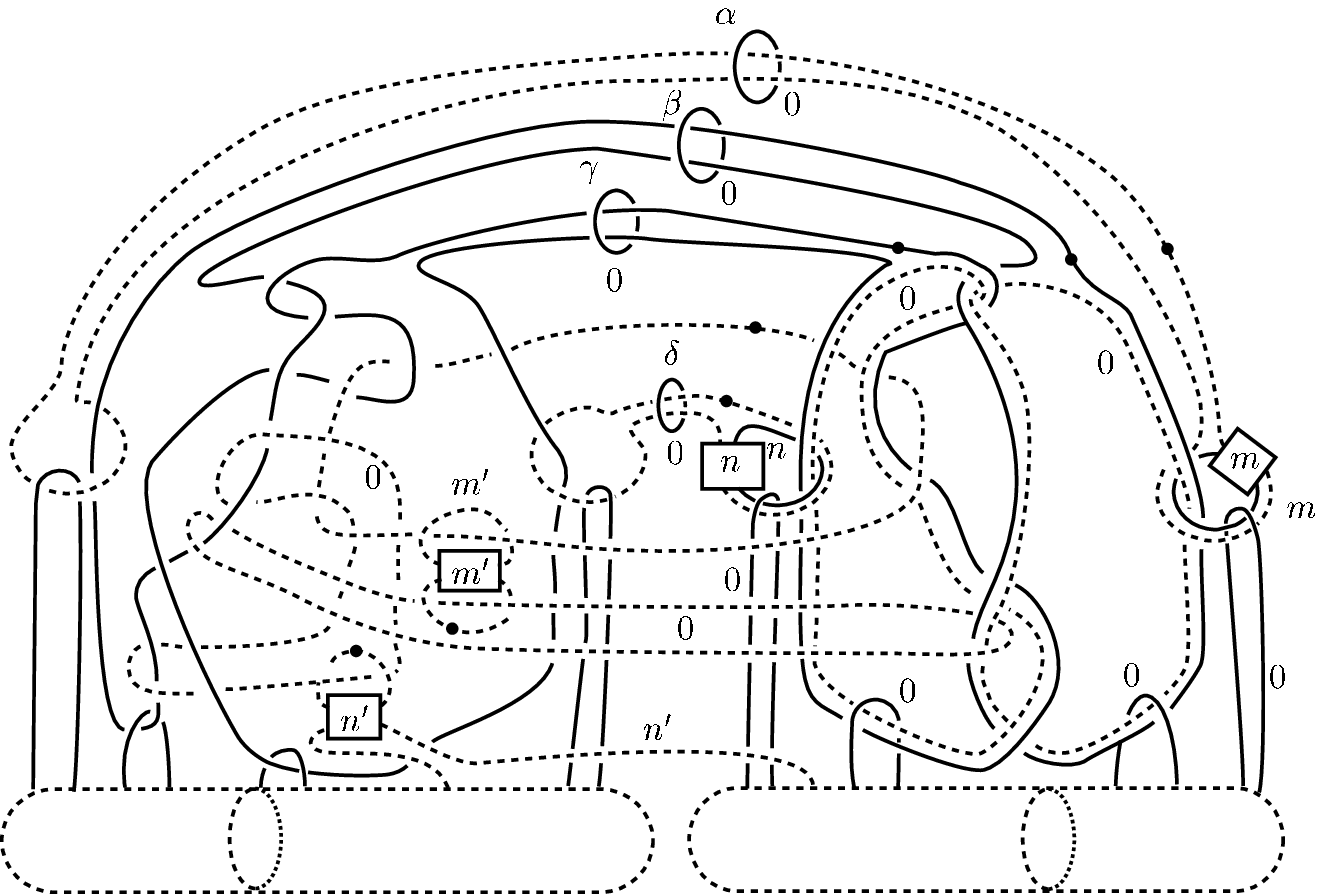}
\caption{${\mathcal S}_{m,n,m',n'}$ canceled.}
\label{can2}
\end{center}\end{figure}

\begin{figure}[htpb]
\begin{center}
\includegraphics{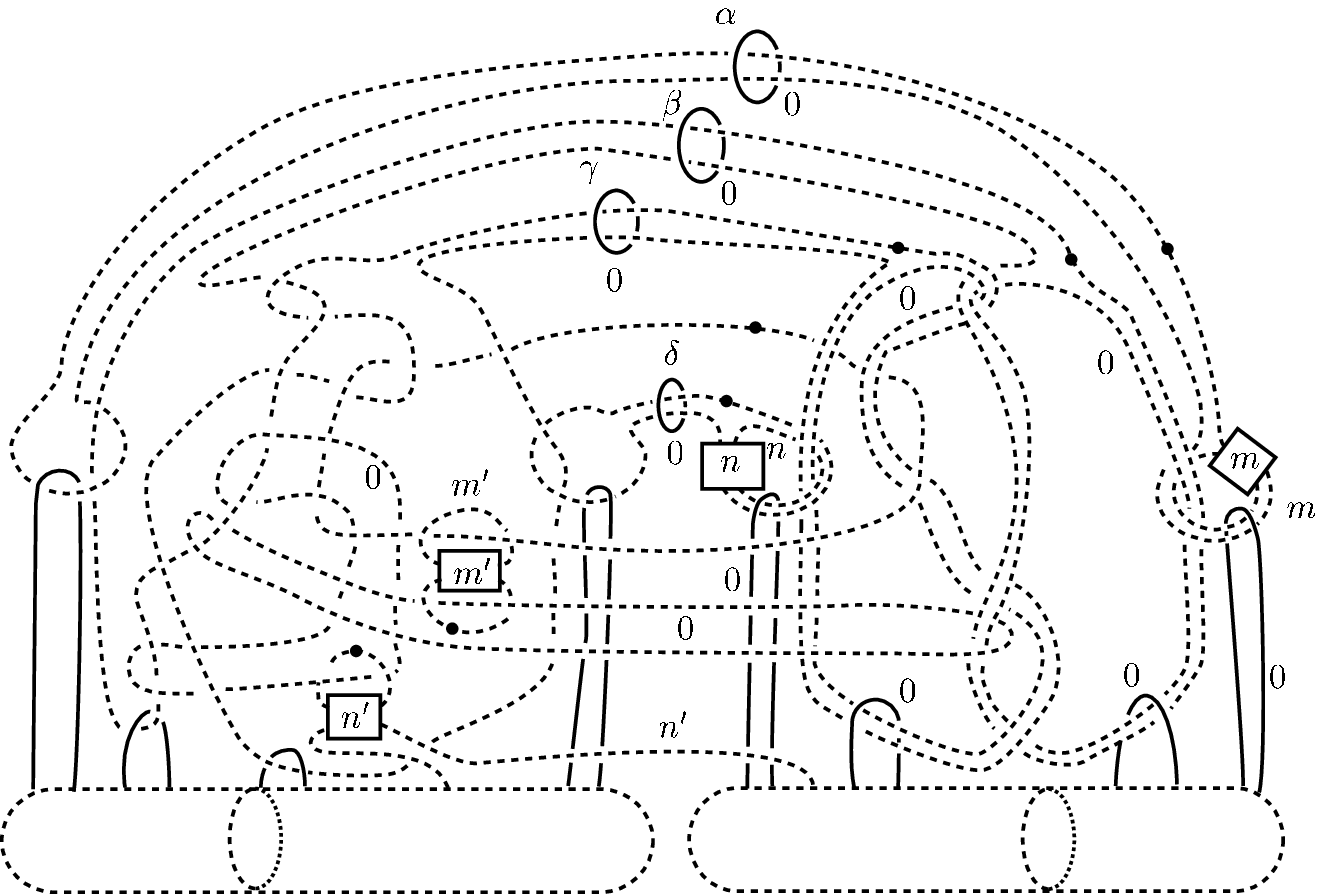}
\caption{${\mathcal S}_{m,n,m',n'}$ canceled.}
\label{can3}
\end{center}\end{figure}

\begin{proposition}
Each of the manifolds ${\mathcal S}_{m,n,m',n'}$ admits a handle decomposition without 1-handles.
In addition the handle decompositon has 4 2-handles.
\end{proposition}
\proof
To prove this lemma, we will find eight 1,2-canceling pairs.
Any canceled pair is drawn by dotted line.
Here the only 1-handle goes on drawing as a ball description.
First we take 4 pairs below as Figure~\ref{can1}.
In addition we take 4 pairs below as Figure~\ref{can2} and \ref{can3}.
Hence we get a handle decomposition 
$${\mathcal S}_{m,n,m',n'}=D^4\cup^8 \text{2-handles} \cup^8\text{3-handles}\cup\text{4-handle}$$

Then, sliding among several canceled handles, we get the figure that 
the attaching circles $\alpha,\beta,\chi,\delta$ are isotopic to the unlink in $\partial D^4$
hence these are canceled out with 4 3-handles in the manifold.

We get a handle decomposition
\begin{equation}
{\mathcal S}_{m,n,m',n'}=D^4\cup^4 \text{2-handles} \cup^4\text{3-handles}\cup\text{4-handle}.
\label{handle}
\end{equation}
\qed

\begin{figure}[htpb]
\begin{center}
\includegraphics{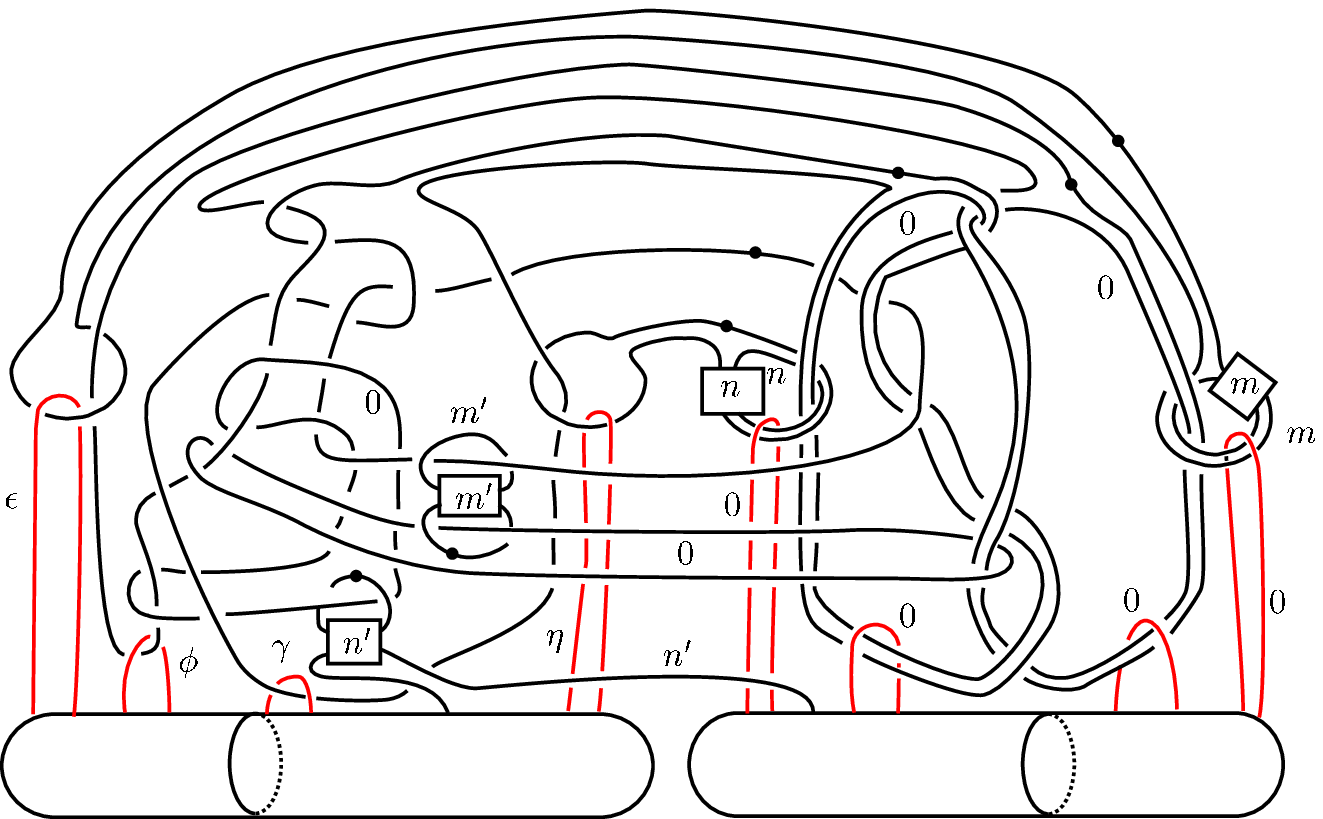}
\caption{${\mathcal S}_{m,n,m',n'}$ canceled.}
\label{step1}
\end{center}
\begin{center}
\includegraphics{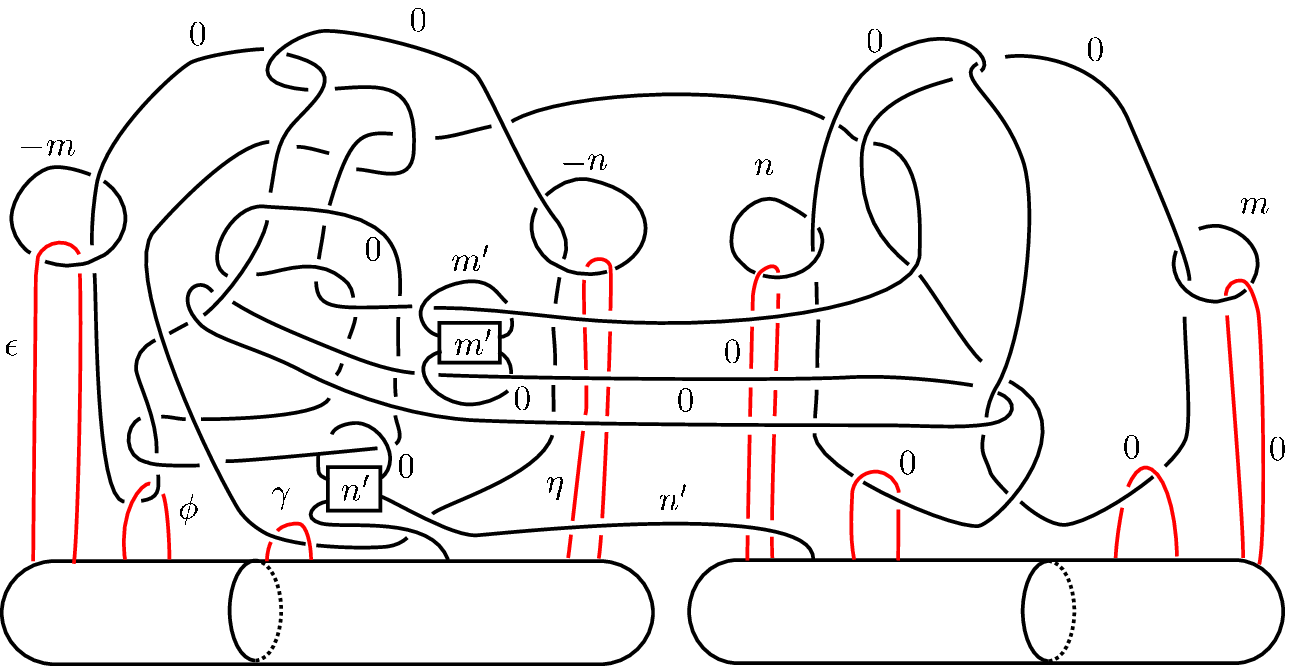}
\caption{2-handles over $\partial D^4$.}
\label{replace}
\end{center}\end{figure}

Thus we get Figure~\ref{step1} as a diagram of ${\mathcal S}_{m,n,m',n'}$.
We put the framed link in $\partial D^4$ of the 2-handles as ${\mathcal F}_{m,n,m',n'}$.

Here four attaching circles (red lines, $\epsilon,\phi,\gamma,\eta$) represent the 2-handles in (\ref{handle}).

\begin{figure}[htpb]
\begin{center}
\includegraphics{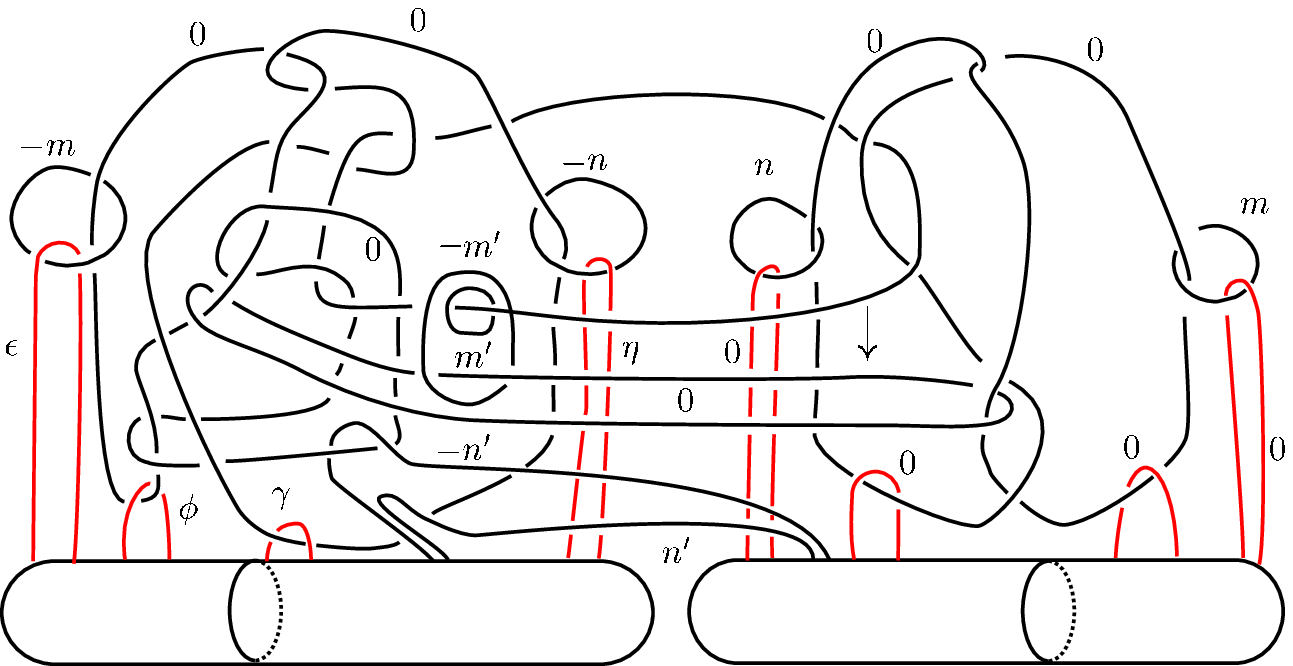}
\caption{2-handles over $\partial D^4$.}
\label{slide1}
\end{center}\end{figure}

\begin{figure}[htpb]
\begin{center}
\includegraphics{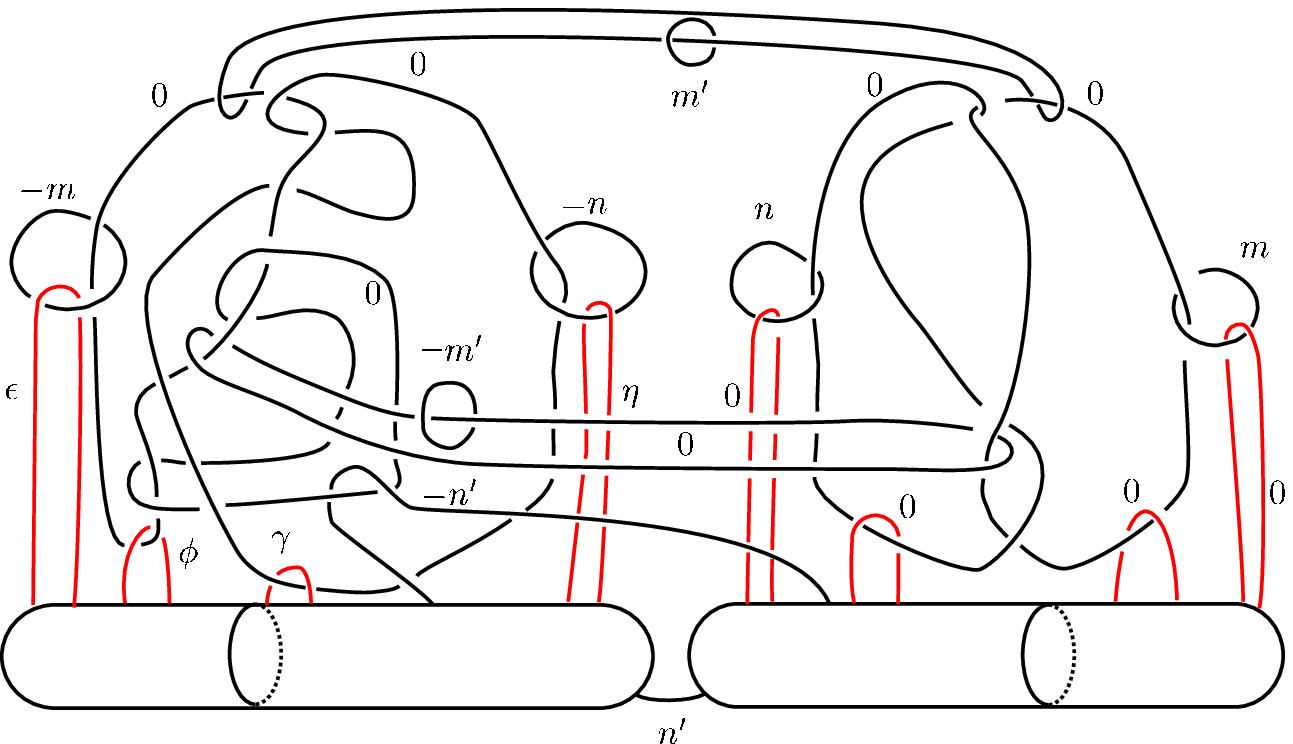}
\caption{2-handles over $\partial D^4$.}
\label{slide2}
\end{center}\end{figure}

\begin{figure}[htpb]
\begin{center}
\includegraphics{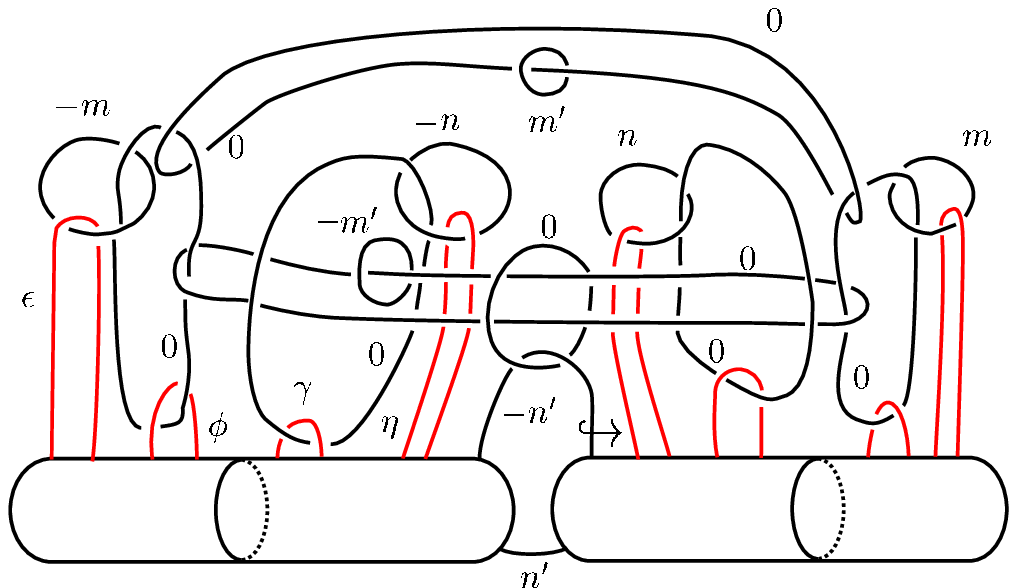}
\caption{2-handles over $\partial D^4$.}
\label{isotopy}
\end{center}\end{figure}

\begin{figure}[htpb]
\begin{center}
\includegraphics{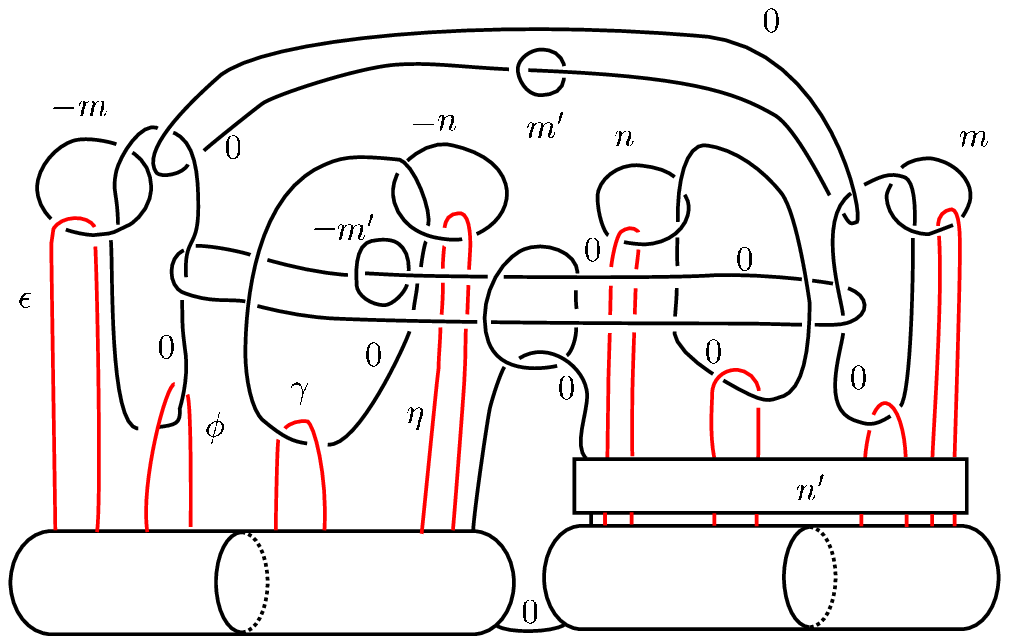}
\caption{2-handles over $\partial D^4$.}
\label{tange}
\end{center}\end{figure}

\begin{figure}[htpb]
\begin{center}
\includegraphics{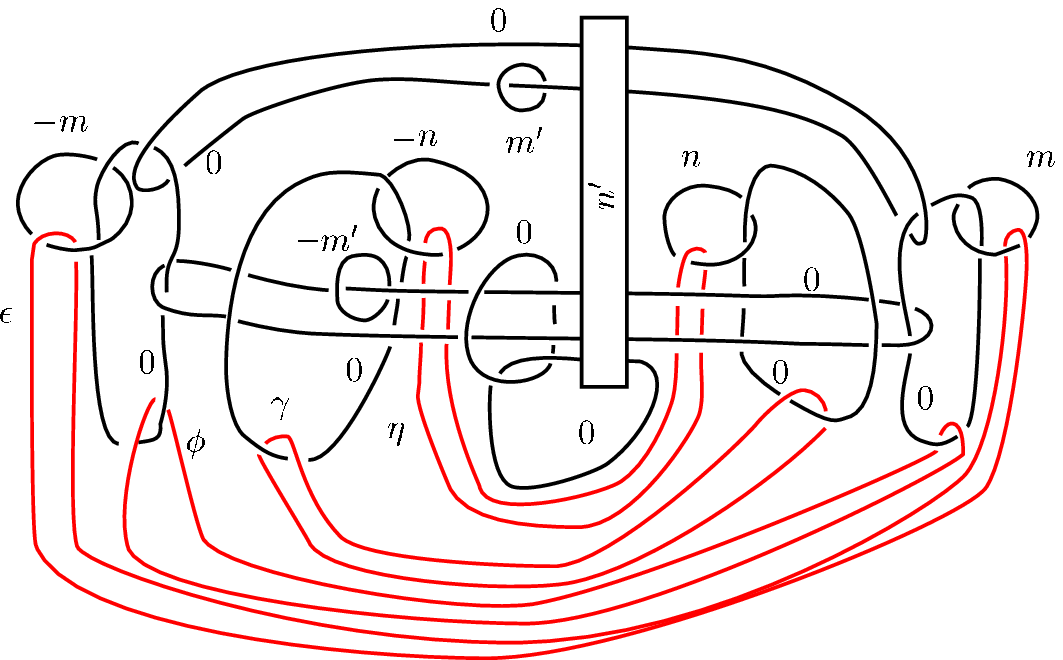}
\caption{2-handles over $\partial D^4$.}
\label{tange2}
\end{center}\end{figure}

\begin{figure}[htpb]
\begin{center}
\includegraphics{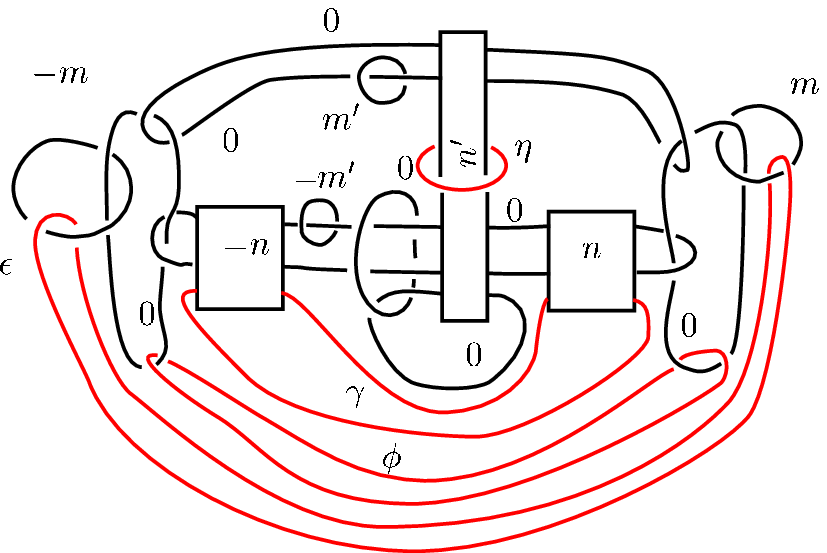}
\caption{2-handles over $\partial D^4$.}
\label{tange3}
\end{center}\end{figure}

\begin{figure}[htpb]
\begin{center}
\includegraphics{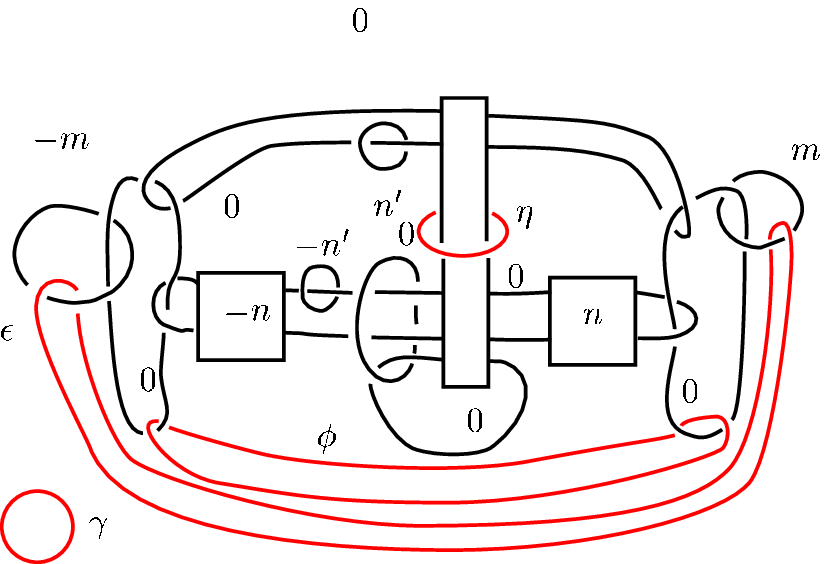}
\caption{2-handles over $\partial D^4$.}
\label{tange4}
\end{center}\end{figure}

\begin{figure}[htpb]
\begin{center}
\includegraphics{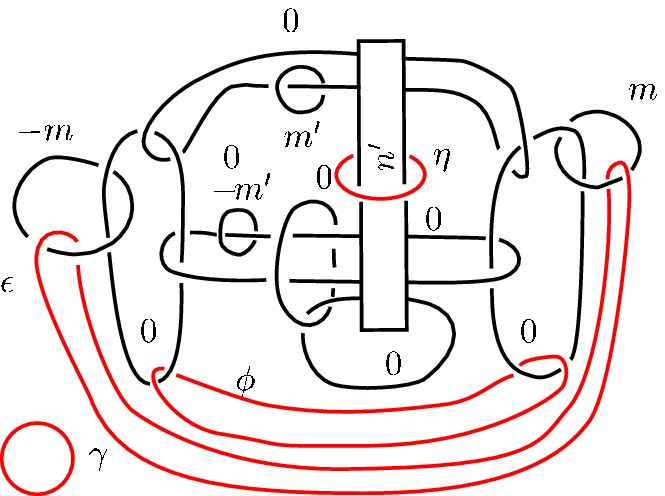}
\caption{2-handles over $\partial D^4$.}
\label{tange5}
\end{center}\end{figure}

\begin{figure}[htpb]
\begin{center}
\includegraphics{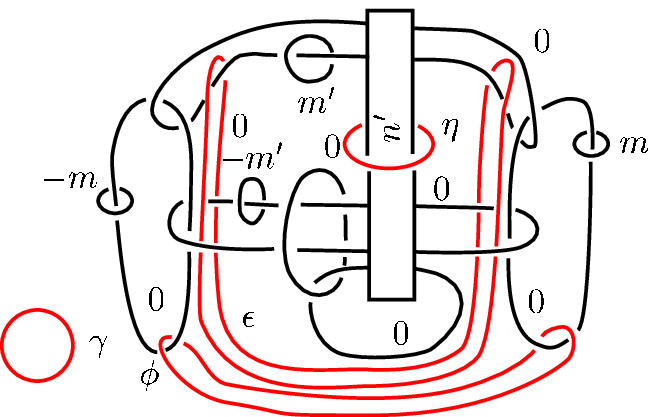}
\caption{2-handles over $\partial D^4$.}
\label{tange6}
\end{center}\end{figure}

\begin{figure}[htpb]
\begin{center}
\includegraphics{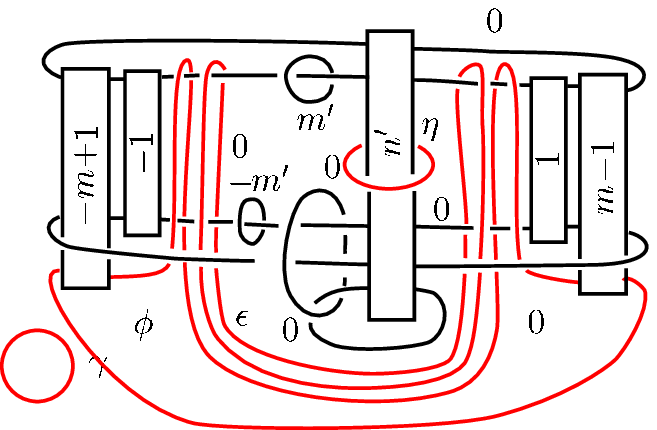}
\caption{2-handles over $\partial D^4$.}
\label{tange7}
\end{center}\end{figure}

\begin{figure}[htpb]
\begin{center}
\includegraphics{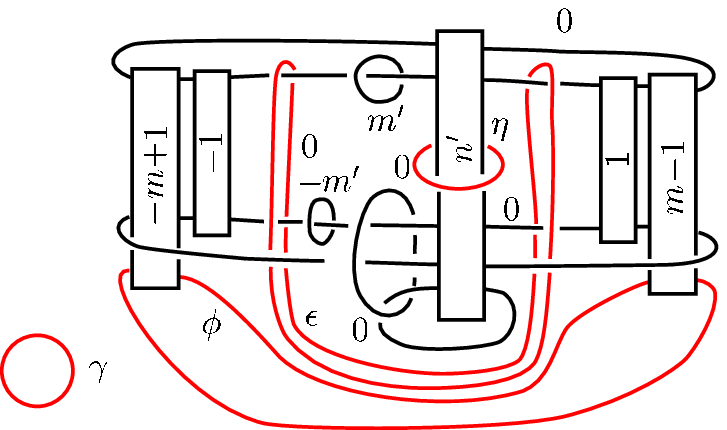}
\caption{2-handles over $\partial D^4$.}
\label{tange8}
\end{center}\end{figure}

Next we show the following.
\begin{lemma}
\label{handleslide}
The 0-framed 2-handles ${\mathcal F}_{m,n,m',n'}$ are, after several handle slides, isotopic to 
${\mathcal F}_{m,0,m',n'}$.
Furthermore two $\gamma,\phi$ of them are separated as 2-component unlink after handle slidings.
\end{lemma}
\proof
Replacing dots of 1-handle to 0-framed 2-handles and sliding handles, we get Figure~\ref{replace}.
Two handle slides give Figure~\ref{slide1}.
Sliding handle as indicated in the figure, we get Figure~\ref{slide2},
and by isotopy we get Figure\ref{isotopy}.
Turning the diagram in the direction of the arrow in Figure~\ref{isotopy} $n'$ times,
we obtain Figure~\ref{tange}.
Removeing bottom canceling 1,2-handle pair, we get Figure~\ref{tange2}.
Sliding central 2-handle, we get Figure~\ref{tange3}.
The curve $\gamma$ in Figure~\ref{tange3} is untied by several handle slidings 
to get a separated 2-handle as in Figure~\ref{tange4}.
At this time the $n$ and $-n$ boxes are untied by rotating (Figure~\ref{tange5}).
Sliding and canceling handles, we get Figure~\ref{tange6} and \ref{tange7}.
In the form we can untie $\epsilon$ by a handle slide as in Figure~\ref{tange8}.
Iterating this process, we get Figure~\ref{tange9}.
\qed

\subsection{Nash's manifolds as a torus surgery.}
\begin{figure}[thbp]
\begin{center}
\includegraphics{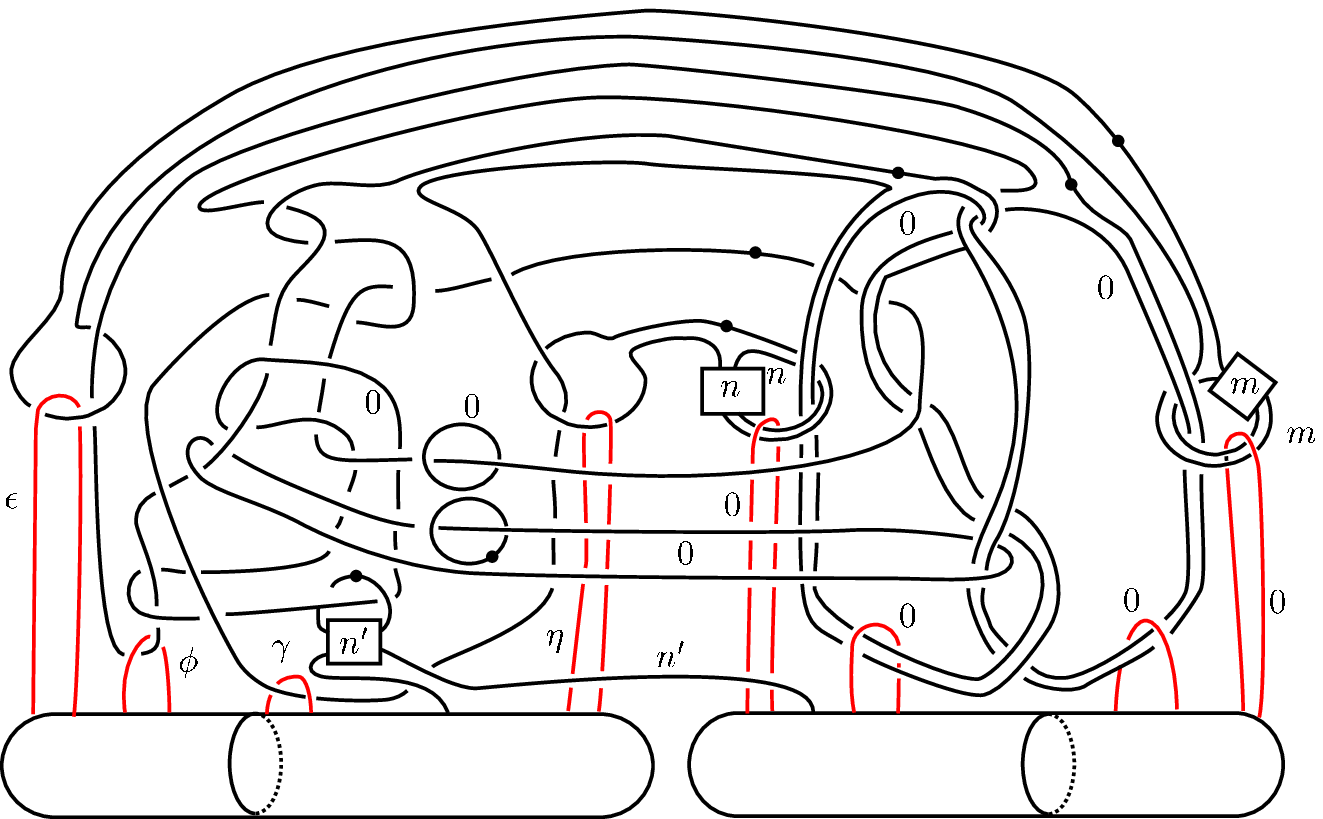}
\caption{${\mathcal S}_{m,n,0,n'}$}
\label{subst0}
\end{center}
\end{figure}

\begin{figure}[thbp]
\begin{center}
\includegraphics{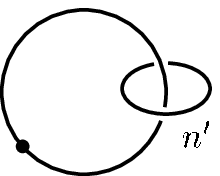}
\caption{${\mathcal S}_{m,n,0,n'}$}
\label{S3S1}
\end{center}
\end{figure}

In the subsection we show that each of Nash's manifold is constructed by a log transformation along a single torus.
\begin{proposition}
\label{zero}
For any $m,n,m',n'$ we have
$${\mathcal S}_{m,n,0,n'}\cong {\mathcal S}_{m,n,m',0}\cong S^4.$$
\end{proposition}
\proof
Putting $m'=0$, we have Figure~\ref{subst0}.
The resulting manifold is the surgering of $S^3\times S^1$ along $\{\text{pt}\}\times S^1$ framing $n'$.
Namely the manifold has the same as Figure~\ref{S3S1}.
This is diffeomorphic to $S^4$.
The manifold ${\mathcal S}_{m,n,m',0}$ is also diffeomorphic to $S^4$ in the similar way.
\qed

As a corollary we have the following.
\begin{corollary}
${\mathcal S}_{m,n,m',n'}$ are given by one log transfomation along a torus.
\end{corollary}
Now we are in a position to prove the main theorem. \\
\subsection{Proof of Theorem~\ref{main}.}
By Lemma~\ref{handleslide} the handle decomposition of ${\mathcal S}_{m,n,m',n'}$ is 
$D^4$ and the same framed link as ${\mathcal F}_{m,0,m',n'}$ and four 3-handles and a 4-handle.
Namely ${\mathcal S}_{m,n,m',n'}$ is the same handle decomposition as ${\mathcal S}_{m,0,m',n'}$.
In particular we have ${\mathcal S}_{m,n,m',n'}\cong {\mathcal S}_{m,0,m',n'}$.
From the Lemma~\ref{exchange} and Proposition~\ref{zero} we have ${\mathcal S}_{m,n,m',n'}\cong S^4$.
\qed

\begin{figure}[htpb]
\begin{center}
\includegraphics{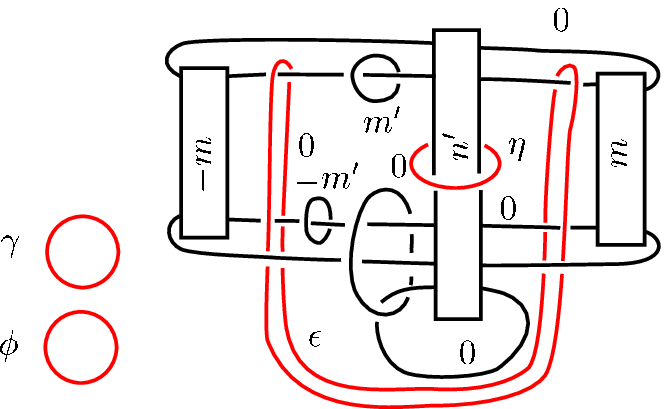}
\caption{2-handles over $\partial D^4$.}
\label{tange9}
\end{center}\end{figure}

\begin{corollary}
\label{cor1}
The diagram Figure~\ref{pcp2} is framed link presentation of $\#^2S^2\times S^1$.
\end{corollary}
\proof
Figure~\ref{tange9} gives a handle decomposition of $S^4$:
$$D^4\cup^2 2\text{-handles}\cup^23\text{-handles}\cup 4\text{-handle}.$$
Therefore the boundary $\partial(D^4\cup^2 2\text{-handles})$ is $\#^2S^2\times S^1$.
\qed\\
This corollary implies $\epsilon$ and $\eta$ in Figure~\ref{pcp2} are
candidates of counterexample of generalized Property R conjecture.

 \noindent
 Motoo Tange\\
 Research Institute for Mathematical Sciences, \\
 Kyoto University, \\
 Kyoto 606-8502, Japan. \\
 tange@kurims.kyoto-u.ac.jp

\begin{figure}[htpb]
\begin{minipage}{.5\textwidth}
\begin{center}
\includegraphics{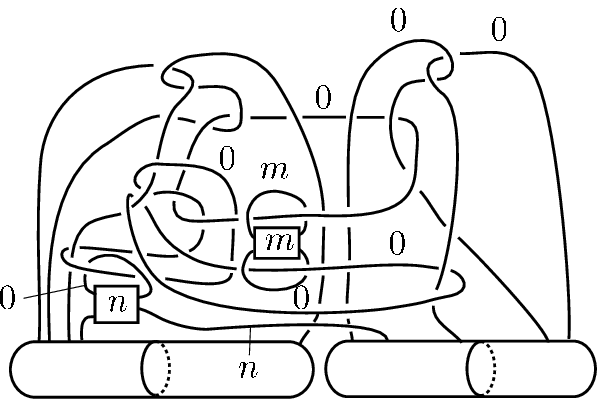}
\caption{}
\end{center}
\end{minipage}
\begin{minipage}{.5\textwidth}
\begin{center}
\includegraphics{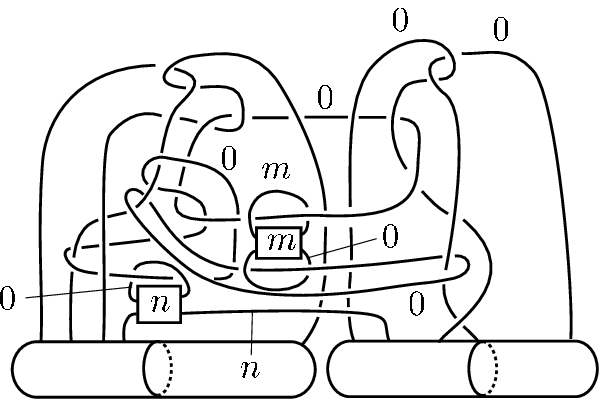}
\caption{}
\end{center}
\end{minipage}
\end{figure}



\begin{figure}[htpb]
\begin{center}
\includegraphics{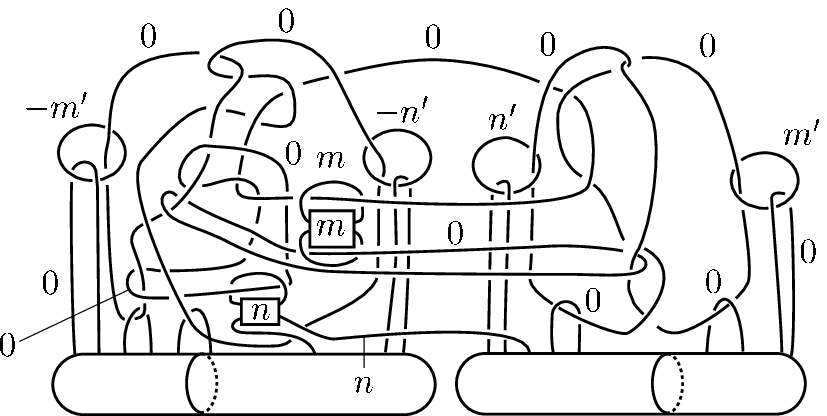}
\caption{$X_{m,n}$}
\label{upd21}
\end{center}\end{figure}


\begin{figure}[htpb]
\begin{center}
\includegraphics{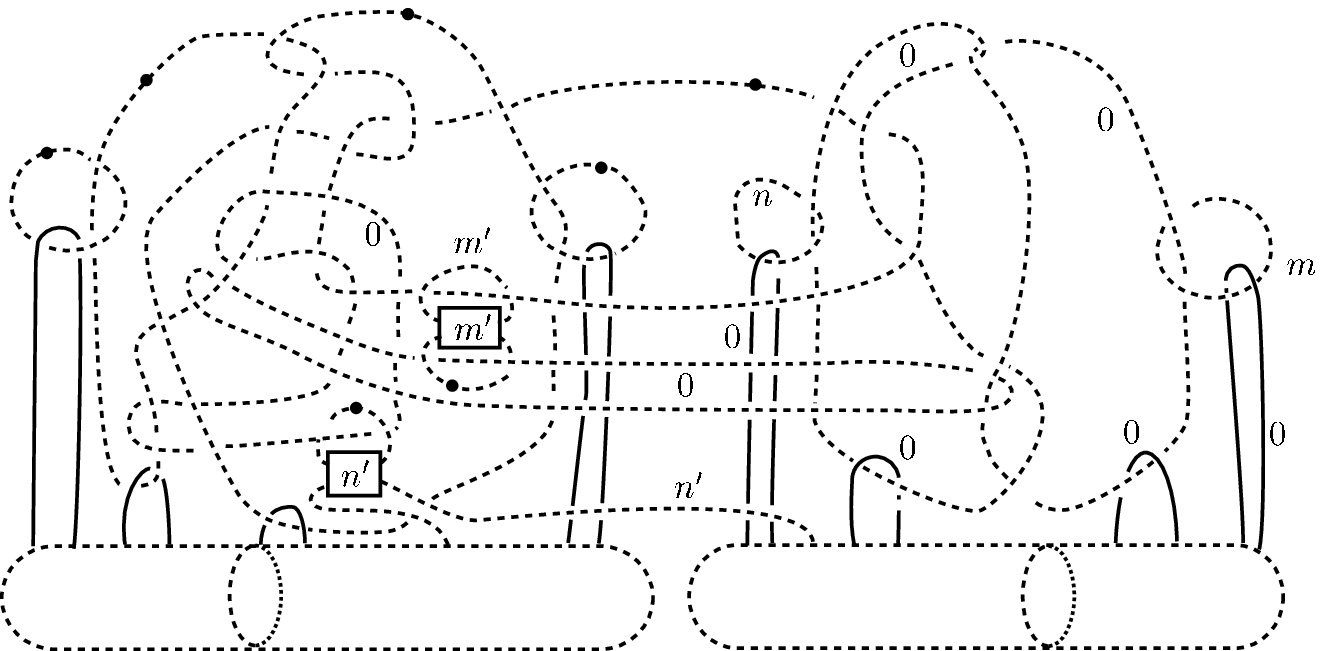}
\caption{${\mathcal S}_{m,n,m',n'}$ canceled.}
\label{can4}
\end{center}\end{figure}
\clearpage




\begin{figure}[thbp]
\begin{center}
\includegraphics{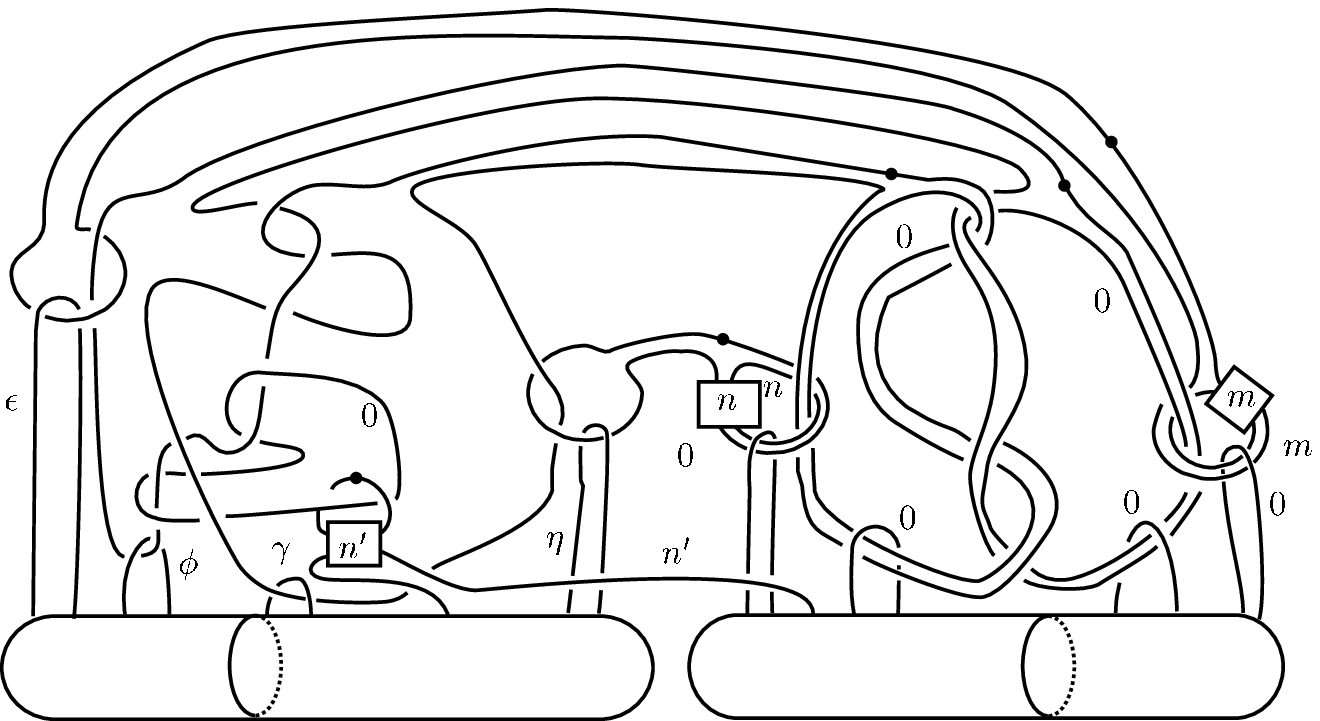}
\caption{${\mathcal S}_{m,n,0,n'}$}
\label{subst1}
\end{center}
\end{figure}

\begin{figure}[thbp]
\begin{center}
\includegraphics{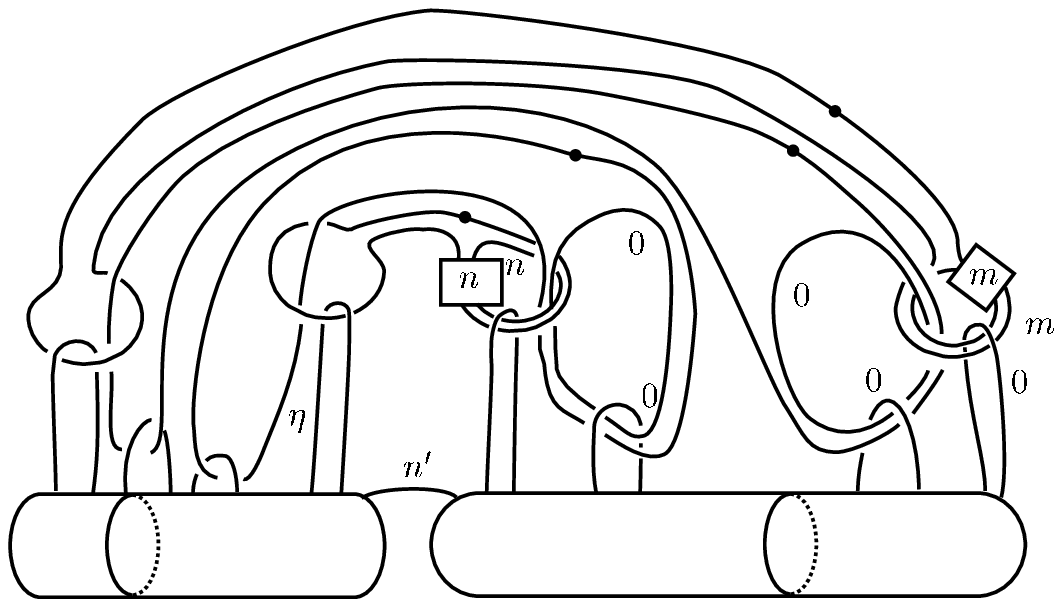}
\caption{${\mathcal S}_{m,n,0,n'}$}
\label{subst2}
\end{center}
\end{figure}

\end{document}